    \input amstex
     \documentstyle{amsppt}
   \pageno=1
\NoRunningHeads

     \magnification = 1200
      \pagewidth{6.5 true in}
    \pageheight{8.7 true in}
   
   \NoBlackBoxes
     \parskip=3 pt
     \parindent = 0.3 true in

\def\({\bigl(}
\def\){\bigr)}
\def\Z{\Bbb Z}
\def\QQ{\Bbb Q}

\def\Q{\Bbb Q}

\def\C{\Bbb C}

\def\A{\Cal A}

\def\N{\operatorname{N}}
\def\Tr{\operatornamewithlimits{Tr}}
\def\Oc{\Cal O}
\def\a{\alpha}

\def\sig{\sigma}

\def\fa{\frak a}

\def\p{\frak p}
\def\P{\frak P}
\def\q{\frak q}
\def\fQ{\frak Q}

\def\Gal{{\operatorname{Gal}}}

\def\coker{{\operatorname{coker}}}
\def\rank{{\operatornamewithlimits{rank}}}

\def\E{{\Cal E}}
\def\S{{\Cal S}}
\def\chibar{{\overline{\chi}}}

\def\Sf{S_{\text{fin}}}

\def\Ta{23}
\def\JT{22}
\def\JTT{21}
\def\Mi{16}
\def\DR{6}
\def\SaBC{19}

\def\Ga{8}
\def\CS{5}
\def\BG{3}
\def\Sn{20}
\def\DST{7}
\def\Co{4}
\def\QuFG{18}
\def\QuFF{17}
\def\MW{15}
\def\Wi{24}
\def\Ko{11}
\def\KoMo{13}
\def\Ka{10}
\def\Ba{1}
\def\Br{2}
\def\JH{9}
\def\Le{14}
\def\Kol{12}
    \topmatter
    \title
Values at $s=-1$ of $L$-functions for multi-quadratic extensions of number fields,
and annihilation of the tame kernel.
    \endtitle
\author
Jonathan W. Sands, University of Vermont  \\
Lloyd D. Simons, Saint Michael's College
\endauthor
\thanks First author partially supported by NSA grant MDA 904-03-1-0003 
\endthanks
\abstract
  Suppose that $\E$ is a totally real number field which is the composite of all of its subfields $E$ that are relative quadratic extensions of a base field $F$. For each such $E$ with ring of integers $\Oc_E$, assume the truth of the Birch-Tate conjecture (which is almost fully established) relating the order of the tame kernel $K_2(\Oc_E)$ to the value of the Dedekind zeta function of $E$ at $s=-1$,
and assume the same for $F$ as well. Excluding a certain rare situation, we prove the annihilation 
of $K_2(\Oc_{\E})$ by a generalized Stickelberger element in the group ring of the Galois group of $\E/F$. Annihilation of the odd part of this group is proved unconditionally.    

\endabstract

    \endtopmatter

     \document

\bigskip\bigskip
           
\heading {I. Notation}
\endheading
 
 Let $F$ be a fixed algebraic number field (finite extension of the rational numbers $\Q$)
and $\E/F$ be a finite abelian extension, with Galois group $G$.    
We fix a set $S$ of primes of $F$ which contains 
all of the infinite primes of $F$, and all of the primes which ramify in $\E$.
For each $\chi$ in the character group $\hat{G}$ of $G$,
we have an Artin 
$L$-function with Euler factors at the
primes in $S$ removed, defined as follows. 
Let $\p$ run through the (finite) primes of $F$ not in $S$,
and $\fa$ run through integral ideals of $F$ which are relatively prime to each of 
the elements of $S$.
Then ${\N} \fa$ denotes the absolute norm of the ideal $\fa$,
 $\sigma_{\fa} \in G$ is the well-defined automorphism 
attached to $\fa$ via the Artin map, and

$${
L^S(s,\chi)=L_{\E/F}^S(s,\chi) =
\sum\Sb {\frak a} \ \text{integral} \\ ({\frak a},S)=1 \endSb
\dfrac{\chi(\sigma_{\fa})}{ {\N} {\frak a}^s}=
\prod\Sb \text{prime } \frak p \notin S \endSb
\bigl(1-\dfrac{\chi(\sigma_{\p})}{ {\N} {\frak p}^{s}}}\bigr)^{-1},$$ 
the Artin $L$-function of the complex variable $s$. Let $\chibar$ be the complex conjugate character of $\chi$.
It is known that $L^S(s,\chi)$ has an analytic continuation to a meromorphic function on $\C$ and a
functional equation relating it to $L^S(1-s,{\chibar})$.

Let $\mu_\infty$ be the group of all roots of 
unity in an algebraic closure of $\E$.
Also let
$${e_{\chi}=\frac{1}{|G|}\sum_{\gamma \in  G}\chi(\gamma)\gamma^{-1}}$$ 
be the primitive idempotent of the group ring $\C[G]$ associated with the character $\chi \in \hat{G}$.
 The Stickelberger function with values 
in $\C[G]$ is defined by 
$${\theta_{\E/F}^S(s)=
\sum_{\chi \in \hat{G}} L^S(s,\chibar)e_{\chi}}.$$ 
The superscript $S$ may be suppressed when $S$ is the minimal set consisting of just the infinite primes and those 
which ramify in $\E/F$. 
Denote by $W_{k+1}(\E)$ the group consisting of all roots of unity which are fixed by the $(k+1)$st power of the
Galois group of $\E(\mu_\infty)/\E$, and let $w_{k+1}(\E)$ be the order of this finite cyclic group. 
By a theorem of Deligne and Ribet [\DR], $w_{k+1}(\E)$$\theta_{\E/F}^S(-k)$ is an element of the integral group ring $\Z[G]$.
Furthermore, Deligne and Ribet show that for each prime $\q$ of $F$ which does not lie in $S$ or contain 
$w_{k+1}(\E)$, the element $$\alpha_{\E/F}^S(\q,k)=(\N\q^{k+1}-\sigma_{\q})\theta_{\E/F}^S(-k)$$ also lies in $\Z[G]$. We will call it an integralized higher Stickelberger element.  

\proclaim{Proposition 1.1 [Sinnott, cf. \Co, Lemma 2.3]} For any integer $C$, the greatest common divisor of the quantities 
$\N(\q)^{k+1}-1$, as $\q$ varies over all
primes of $F$ not dividing $w_{k+1}(F)C$, is $w_{k+1}(F)$.
\endproclaim

\heading{II. Conjectures on Annihilation of $K$-groups}
\endheading

For any ring $R$, $K_n(R)$ will denote the $n$th algebraic $K$-group of $R$, as defined by Quillen.
Our focus will be on the ring $\Oc_\E$ of algebraic integers of $\E$, and the ring of $S$-integers $\Oc_\E^S$, which consists of the elements of $\E$ whose valuation is non-negative at every prime not above one of those in $S$. 
It is known (see for example Corollary 1.11 of [\Mi]) that $K_0(\Oc_\E^S)$ is isomorphic to the direct sum of $\Z$ and the class group of $\Oc_\E^S$.
Using this isomorphism, a key prediction of the Brumer-Stark conjecture [\Ta] states that the group ring element 
$\alpha_{\E/F}^S(\q,0)$ annihilates the Galois module $K_0(\Oc_\E^S)$ when $w_1(\E)\notin \q \notin S$, and $S$ has cardinality greater than 1. The conjecture follows from Stickelberger's theorem
when $F=\Q$, and is proved for $\E/F$ relative quadratic in [\Ta]. With minor exceptions, the conjecture is proved for
$\E/F$ multiquadratic in [\DST].

Coates and Sinnott prove in [\CS] that 
$\alpha_{\E/F}^S(\q,1)$ annihilates the odd part of $K_2(\Oc_{\E})$ when 
$F=\Q$, $w_2(\E)\notin \q \notin S$, and $S$ is minimal. 
In 1977, Coates [\Co] suggested that the natural generalization of the above phenomena would be that  
$\alpha_{\E/F}^S(\q,k)$ annihilates $K_{2k}(\Oc_{\E}^S)$ when $F=\Q$ and $w_{k+1}(\E)\notin \q \notin S$.
As an example of the evidence which has accumulated since then, Banaszak proved in [\Ba] that for $F=\Q$ and arbitrary non-negative $k$, $\alpha_{\E/F}^S(\q,k)$ annihilates the odd part of the subgroup of $K_{2k}(\Oc_{\E})$ consisting of elements whose image in $K_{2k}(\E)$ is divisible. More recent results such as those in [\SaBC], [\Sn], [\BG] have provided evidence for general $F$.
For example, [\SaBC] shows that the results of [\CS]   and  [\Ba]  mentioned above for $F=\Q$ also hold after a base change: for the same $\E$ and $S$, the base field $F$ may be taken as any intermediate field between $\Q$ and $\E$. 

The well-known conjecture of Birch and Tate (see section 4 of [\JTT]) may be viewed as strengthening the special case of the above when
$\E=F$ and $S$ is minimal. In this case, $\chi$ must be the trivial character $\chi_0$ and
we have $L_{F/F}(s,\chi_0)=\zeta_F(s)$, the Dedekind zeta function of $F$. It is known from [\Ga] and [\QuFG] that $K_2(\Oc_F)$ is finite.

\proclaim{Conjecture(Birch-Tate)}
Suppose that $F$ is totally real of degree $r_1(F)$. Then 
$$\zeta_F(-1)= (-1)^{r_1(F)} \dfrac{|K_2(\Oc_F)|}{w_2(F)}$$
\endproclaim

 Deep results on Iwasawa's Main conjecture in [\MW] and [\Wi] lead to the following
(see [\Ko]).

\proclaim{Theorem 2.1}
The Birch-Tate Conjecture holds if $F$ is abelian over $\Q$, and the odd part
holds for all $F$.
\endproclaim 

This result will provide the key connection between $L$-function values and 
$K$-groups in the proof of our main result.

\heading{III.} The Structure of $K_2(\Oc_{\E}^S)$.
\endheading

 Recall that $\E/F$ is an abelian extension and $S$ is a finite set of primes of $F$, containing all of the infinite primes of $F$ and the primes which ramify in $\E$. Then $S_\E$ will denote primes of $\E$ which lie over those in $F$, and $\Oc_\E^S=\Oc_\E^{S_\E}$. We now review how the group $K_2(\Oc_{\E}^S)$ may be viewed as a subgroup of $K_2(\E)$, and how Matsumoto's theorem
gives a useful description of $K_2(\E)$ which actually holds for any field.

\proclaim{Theorem 3.1 (Matsumoto)}
$K_2(\E)\cong \E^\times
\otimes_\Z \E^\times /\langle\{x\otimes(1-x):1\neq x \in \E^\times\}\rangle$.
\endproclaim

 $K_2(\E)$ is therefore generated by Steinberg symbols $\{\alpha,\beta\}$ for $\alpha, \beta \in  \E^\times$;
such an element corresponds to the coset of $\alpha\otimes\beta$ under the isomorphism of the theorem.

For each prime ideal $\P$ of $\Oc_\E$, let $\nu_\P$ be the corresponding discrete valuation on $\E$. Then for $\{\alpha,\beta\}\in K_2(\E)$, the tame symbol $(\alpha,\beta)_\P$ is a 
well-defined element in $(\Oc_\E/\P)^\times$ given by 
$$(\alpha,\beta)_\P=(-1)^{\nu_\P(\alpha)\nu_\P(\beta)}\alpha^{\nu_\P(\beta)}/\beta^{\nu_\P(\alpha)}\pmod{\P}.$$

The following well-known consequence of work of Bass--Tate and
Quillen (see Quillen, [\QuFF]) gives a useful characterization of $K_2(\Oc_\E^\S)$.

\proclaim{Theorem 3.2 (Localization Sequence)}
There is a natural exact sequence of abelian groups
$$ 0 \rightarrow K_2(\Oc_\E^S)\rightarrow K_2(\E)\rightarrow \oplus_{\P\notin S_\E}(\Oc_\E/\P)^\times 
\rightarrow 0$$
 where the homomorphism with domain
$K_2(\E)$ sends $\{\alpha,\beta\}$ to the element whose component at $\P$ is the tame symbol 
$(\alpha,\beta)_\P$.
\endproclaim

Setting $\zeta_F^S(s)=L_{F/F}^S(s,\chi_0)$, and $\Sf$ equal to the set of finite primes in $S$, we observe an immediate corollary which is perfectly analogous to the situation for $k=0$.

\proclaim{Corollary 3.3 (Extended Birch-Tate Conjecture)}
If the Birch-Tate conjecture holds for $F$, then for any $S$ containing all of the infinite primes, we have 
$$\zeta_F^S(-1)=(-1)^{|S|} \dfrac{|K_2(\Oc_F^S)|}{w_2(F)}.$$ Equality holds unconditionally for the odd part.
\endproclaim
\demo{Proof} Comparing the localization sequences for $\Oc_F$ and $\Oc_F^S$ (by direct argument, or by naturality and the Snake lemma) leads to the short exact sequence
$$ 0 \rightarrow K_2(\Oc_F)\rightarrow K_2(\Oc_F^S)\rightarrow \oplus_{\p\in\Sf}(\Oc_F/\p)^\times 
\rightarrow 0.$$
From this we see that $$|K_2(\Oc_F^S)|=|K_2(\Oc_F)|\prod_{\p\in\Sf}(\N\p-1).$$
At the same time, the Euler product formula shows that $$\zeta_F^S(-1)=\zeta_F(-1)\prod_{\p\in\Sf}(1-\N\p).$$
Thus introducing $S$ contributes factors of the same absolute value to both sides of the equation in the original Birch-Tate conjecture, and also contibutes a change of sign for each finite prime in $S$. The result follows, using Theorem 2.1 for the odd part. 
\qed
\enddemo
 
 We will also need to use a well-known corollary (stated as Corollary 3.6 below) of a result of Kahn [\Ka].
For completeness and consistency of exposition, we sketch the proof here. See [\Kol, Remark 2.9] for a more general version of
this result involving \'etale cohomology. 
  
For a $G$-module $M$, we denote the quotient module of co-invariants by $M_G$. 

\proclaim{Theorem 3.4 (Kahn [\Ka, 5.1])} Let $\E/F$ be a Galois extension of number fields with Galois group $G$. Let
$r_{\infty}(\E/F)$ be the number of real places of $F$ which extend to complex places of $\E$, and consider $K_2(\E)$ as a $G$-module for the natural action of $G$. Then the transfer map from $K_2(\E)$ to $K_2(F)$ induces 
a map which fits into an exact sequence:
$$0\rightarrow K_2(\E)_G\rightarrow K_2(F)\rightarrow \{\pm 1\}^{r_\infty(\E/F)}\rightarrow 0.$$
\endproclaim
 
We first prove a lemma.
\proclaim{Lemma 3.5} With the assumptions as in the theorem, also let $S$ be a finite set of primes of $F$ containing the infinite primes and those which ramify from 
$F$ to $\E$. Then there is a natural exact sequence  
$$0\rightarrow K_2(\Oc_\E^S)_G\rightarrow K_2(\E)_G\rightarrow \oplus_{\p\notin S}(\Oc_F/\p)^\times \rightarrow 0.$$
\endproclaim
\demo{Proof} Begin with the exact localization sequence
$$ 0 \rightarrow K_2(\Oc_\E^S)\rightarrow K_2(\E)\rightarrow \oplus_{\P\notin S_\E}(\Oc_\E/\P)^\times 
\rightarrow 0.$$
This is an exact sequence of $G$-modules, given the natural action of $G$ on each term. The corresponding long exact sequence in homology becomes
$$H_1(G, \oplus_{\P\notin\S_\E}(\Oc_\E/\P)^\times)
\rightarrow K_2(\Oc_\E^S)_G\rightarrow K_2(\E)_G\rightarrow \bigl(\oplus_{\P\notin S_\E}(\Oc_\E/\P)^\times\bigr)_G 
\rightarrow 0.$$
 We have $\oplus_{\P\notin S_\E}(\Oc_\E/\P)^\times=
\oplus_{\p\notin S}\oplus_{\P | \p}(\Oc_\E/\P)^\times$,
and $\oplus_{\P | \p}(\Oc_\E/\P)^\times$ is a $G$-module
induced from the module $(\Oc_\E/\P)^\times$ over the
decomposition group $G_\P$ of a single fixed $\P$. By Shapiro's lemma, we may then compute the homology of this summand over $G_\P$,
and $G_\P$ can be identified with the Galois group of $\Oc_\E/\P$ over $\Oc_F/\p$, since $\P$ is unramified over $\p$. But the Tate cohomology of 
multiplicative groups of finite fields is trivial. This implies that $H_1(G, \oplus_{\P\notin\S_\E}(\Oc_\E/\P)^\times)=0$
and 
$\bigl(\oplus_{\P\notin\S_\E}(\Oc_\E/\P)^\times\bigr)_G$
is isomorphic to the image of the norm map, which is $\oplus_{\p\notin S}(\Oc_F/\p)^\times.$ 
\qed 
\enddemo 
\proclaim{Corollary 3.6 ([\KoMo, Prop. 1.6])} With the assumptions as in the lemma, the transfer map 
from $K_2(\Oc_\E^S)$ to $K_2(\Oc_F^S)$ induces a homomorphism which fits into a natural exact sequence
$$0 \rightarrow K_2(\Oc_\E^S)_G\rightarrow K_2(\Oc_F^S)\rightarrow \{\pm 1\}^{r_\infty(\E/F)}\rightarrow 0 $$ 
\endproclaim
\demo{Proof}
 Applying the lemma to both $\E/F$ and $F/F$ and using naturality yields the commutative diagram
$$
\CD
0 @>>> K_2(\Oc_\E^S)_G @>>> K_2(\E)_G @>>> \oplus_{\p\notin S}(\Oc/\p)^\times @>>> 0 \\ 
  @.   @V{Tr}VV              @V{Tr}VV        @VVV                          \\
0 @>>> K_2(\Oc_F^S) @>>> K_2(F) @>>> \oplus_{\p\notin S}(\Oc/\p)^\times @>>> 0 
\endCD
 $$
One checks (again using the fact that each prime $\p\notin S$ is unramified in $\E/F$) that the
 vertical map on the right is an isomorphism. The snake lemma then allows us to identify the kernel and cokernel of the vertical map on the left with the kernel and cokernel of the vertical map in the middle. By Kahn's theorem, these are
$0$ and $\{\pm 1\}^{r_\infty(\E/F)}$, respectively. The result follows.
\qed
\enddemo 

\heading{IV. Statements of Results}
\endheading

From now on, we assume that $\E$ is a totally real multiquadratic extension of $F$, that is, a composite of a finite number of totally real relative quadratic extensions $E$ of the number field $F$. We further assume that $S$ is a finite set of primes of $F$ which contains all of the infinite primes as well as all of the primes that ramify in $\E/F$. Suppose also that $\q$ is a fixed prime of $F$ such that $w_{2}(\E)\notin \q \notin S$.
Let $S_2$ denote the set of primes of
$F$ that lie above $2$ (the ``dyadic primes'' of $F$).
Let $G$ be the Galois group of $\E/F$, an elementary abelian 2-group of rank $m=\rank_2(G)$. 

Our main results are the following.
\proclaim{Theorem 4.1}
 Let $\E$ be any totally real multiquadratic extension of $F$. Let $S$ and $\q$ be as just described. 
\item{a.} Then $\alpha_{\E/F}^S(\q,1)$ annihilates the odd part of $K_{2}(\Oc_{\E}^S)$.
\item{b.} If the full Birch-Tate conjecture holds for $F$ and for each relative quadratic extension $E$ of $F$ inside $\E$, then 
$2\alpha_{\E/F}^S(\q,1)$ annihilates $K_{2}(\Oc_{\E}^S)$. 
\item{c.} If $\E(\sqrt{-1})$ is not the maximal multiquadratic extension of $F$ which is unramified outside of $S\cup S_2$,
and the Birch-Tate conjecture holds for $F$ and for each relative quadratic extension $E$ of $F$ inside $\E$, then 
$\alpha_{\E/F}^S(\q,1)$ annihilates $K_{2}(\Oc_{\E}^S)$. 
\endproclaim
  
 From 4.1 and 2.1 we obtain the following corollary. For multiquadratic extensions, it strengthens the result that would follow from [\CS] and [\SaBC] by allowing smaller sets $S$ and by encompassing the 2-part. 
 
\proclaim{Corollary 4.2} Suppose $\E$ is abelian over $\Q$ and $\E(\sqrt{-1})$ is not the maximal multiquadratic extension of $F$ which is unramified outside of $S\cup S_2$. Then
$\alpha_{\E/F}^S(\q,1)$ annihilates $K_{2}(\Oc_{\E}^S)$.
\endproclaim
 
\proclaim{Remark 4.3} {\rm If $F$ does not possess a system of totally positive fundamental $S\cup S_2$-units, in particular, if $F$ does not possess a system of totally positive fundamental units, then 
the condition on $\E(\sqrt{-1})$ is automatically met for any $\E$. For if $\eta$ is a fundamental $S\cup S_2$-unit of $F$ 
such that neither $\pm \eta$ is totally positive, then $\E(\sqrt{-1},\sqrt{\eta})$ is a multiquadratic extension of $F$ which is unramified outside of $S\cup S_2$ and strictly larger than $\E(\sqrt{-1})$, as it contains the extension $\E(\sqrt{\eta})$ which is neither totally real nor totally complex. For example, if $F$ contains $\sqrt{2}$, then the condition is met because the unit $\eta=\sqrt{2}-1$ lies in $F$. Even among totally real fields of degree 5, systems of totally positive units are rare. A recent extensive computer search [\JH] found the first examples of such fields, the smallest discriminant of those found so far being $405,673,292,473$, and others occurring about once in every ten million totally real quintic fields considered.} 
\endproclaim

\proclaim{Remark 4.4} {\rm Theorem 4.1 implies that, given the stated hypotheses, $\alpha_{\E/F}^S(\q,1)$ annihilates $K_{2}(\Oc_{\E})$, for (as seen in the proof of Corollary 3.3) $K_{2}(\Oc_{\E})$ may be identified with a subgroup of $K_{2}(\Oc_{\E}^S)$. Compare this with the very different situation of the Brumer-Stark conjecture in which $K_{0}(\Oc_{\E}^S)$ is a quotient of $K_{0}(\Oc_{\E})$.}
\endproclaim

\heading{V. A formula for $\alpha_{\E/F}^S(\q,1)$}
\endheading
  We begin the proof of Theorem 4.1 by obtaining an expression for  $\alpha_{\E/F}^S(\q,1)$ in terms of 
the invariants of the relative quadratic extensions $E/F$. Let $\chi_0$ be the principal character of the Galois group
$G$ of $\E/F$ (of order $2^m$). Put $T_{\E/F}=\sum_{\sigma\in G} \sig = 2^m e_{\chi_0}$, and for each non-principal character $\chi$, let $E_\chi$ be the fixed field of the kernel $\ker(\chi)$ of $\chi$ (so $E_\chi/F$ is a relative quadratic extension). We put
$T_{\E/E_\chi}=\sum_{\sigma\in \ker(\chi)} \sigma$ and fix a choice of $\tau_\chi\in G$ such that $\tau_\chi$ restricts to the nontrivial element of the Galois group of $\E_\chi/F$. Then $e_\chi=2^{-m}(1-\tau_\chi)T_{\E/E_\chi} $. 
Thus  in the group ring ${\QQ}[G]$ we have the equality
$$\theta_{\E/F}^S(-1)=2^{-m}L^S(-1,\chi_0)T_{\E/F}
+2^{-m}\sum_{\chi\neq \chi_0}L^S(-1,\chi)(1-\tau_{\chi}) T_{\E/E_{\chi}}.$$ 
Now
$L^S(-1,\chi_0)=\zeta_F^S(-1)$, and $L^S(-1,\chi)=\zeta_{E_\chi}^S(-1)/\zeta_{F}^S(-1)$, by the standard properties of
Artin $L$-functions. Assuming the Birch-Tate conjecture for $F$ and each $E_{\chi}$, it follows from Corollary 3.3 that
$$(-1)^{|S|}\theta_{\E/F}^S(-1)=2^{-m}\dfrac{|K_2(\Oc_F^S)|}{w_2(F)} T_{\E/F}
+2^{-m}\sum_{\chi\neq \chi_0} (-1)^{|S_{E_{\chi}}|}\dfrac{w_2(F)}{w_2(E\chi)}\dfrac{|K_2(\Oc^S_{E_\chi})|}
{|K_2(\Oc^S_F)|}(1-\tau_{\chi}) T_{\E/E_{\chi}},$$
and indeed equality of the odd parts holds unconditionally. 
To obtain $\alpha_{\E/F}^S(\q,1)$, we now just multiply by
$\N\q^{2}-\sigma_\q$, noting that $\sigma_\q T_{\E/F}=T_{\E/F}$.
 
\proclaim{Proposition 5.1} 
Suppose that $w_2(\E)\notin \q \notin S$.
Assuming the Birch-Tate conjecture holds for each $E_\chi$ and for $F$, we have
$$
\multline
(-1)^{|S|}\alpha_{\E/F}^S(\q,1)= \\
\dfrac{(\N\q^2-1)}{w_2(F)}\dfrac{|K_2(\Oc_F^S)|}{2^m}T_{\E/F}
+ \sum_{\chi\neq \chi_0}(-1)^{|S_{E_\chi}|} \dfrac{w_2(F)}{w_2(E_\chi)}\dfrac{|K_2(\Oc^S_{E_\chi})|}
{2^m|K_2(\Oc^S_F)|}(\N\q^2-\sigma_\q)(1-\tau_{\chi}) T_{\E/E_{\chi}}.
\endmultline$$
Equality of the odd parts holds unconditionally.
\endproclaim

 Our proof of Theorem 4.1 will consist in showing that, after multiplication by the appropriate power of 2 for each part of the theorem, each individual term of the sum in Proposition 5.1 is integral and annihilates $K_2(\Oc_{\E}^S)$. 

\heading{VI. Orders and 2-ranks of K-groups}
\endheading

We proceed by computing the orders and 2-ranks of 
the relevant $K$-groups via Kummer Theory.
Let $\Gamma_F^S=\{a\in F^\times \, : \nu_\p(a) \text{ is even for every prime } \p \notin S\cup S_2\}$. For an abelian group
$A$, we let $_2A$ denote the subgroup of $A$ consisting of those elements whose order divides 2.

 The following Proposition may be derived from a Bockstein sequence. It suits our purposes better to give a fairly explicit proof.
\proclaim{Proposition 6.1} The group $_2K_2(\Oc_F^S)$ is isomorphic to a quotient of $\Gamma_F^S/(F^\times)^2$ 
by a subgroup $D_F/(F^\times)^2$ of order 2. This subgroup is generated by the coset of $\sqrt{b}$, where $F(\sqrt{b}$ 
is the first layer of the cyclotomic $\Z_2$-extension of $F$.
\endproclaim

\demo{Proof} (see [\Br, Theorem 3]) A special case of Tate's Theorem 6.1 in [\JT] implies that there is a surjective homomorphism 
$\varphi: F^\times \rightarrow  {}_2K_2(F)$ defined by $\varphi(a)=\{-1,a\}$. Considering the localization sequence,
$$0\rightarrow K_2(\Oc_F^S) \rightarrow K_2(F)\rightarrow \oplus_{\p\notin S}(\Oc_F/\p)^\times \rightarrow 0,$$
we see that $\{-1,a\}$ is in the kernel of the homomorphism with domain $K_2(F)$ if and only if $(-1,a)_\p=(-1)^{\nu_\p(a)}$ equals 1 
for each $\p\notin S\cup S_2$. (For $\p\in S_2$, $(-1,a)_\p=1$ because $-1\equiv 1 \pmod{\p}$). Thus the elements of order
dividing 2
in the kernel are given by exactly $\{\{-1,a\} : a \in \Gamma_F^S\}$. The exactness of the sequence allows us to 
identify this kernel with $_2K_2(\Oc_F^S)$. Thus $\varphi$ induces a surjective homomorphism from $\Gamma_F^S/(F^\times)^2$ to
$_2K_2(\Oc_F^S)$. Tate's Theorem 6.3 of [\JT] computes the index of 
$(F^\times)^2$ in the kernel $D_F$ of $\varphi$ on $F^\times$ as $2^{r_2(F)+1}$, 
where $r_2(F)$ is the number of pairs of complex embeddings of $F$. In our situation, $r_2(F)=0$, and 
(see Lescop [\Le, Prop. 2.4]) a generator of $D_F/(F^\times)^2$ is
represented by $b$ such that $F(\sqrt{b})$ is the first layer
$F_2$ of the cyclotomic $\Z_2$-extension of $F$. Since this is a quadratic extension unramified outside of $S_2$, we have
$b\in \Gamma_F^S$. Thus $D_F\subset \Gamma_F^S$ and 
$_2K_2(\Oc_F^S)$ is isomorphic to the quotient of $\Gamma_F^S/(F^\times)^2$ by the subgroup $D_F/(F^\times)^2$ of order 2.
\qed
\enddemo
 
 \proclaim{Lemma 6.2} The group $\Gamma_F^S/(F^\times)^2$ has exponent 2 and order divisible by $2^{m+1}$.
If $\E(\sqrt{-1})$ is not the maximal multiquadratic extension of $F$ which is unramified outside of 
$S\cup S_2$, then  $\Gamma_F^S/(F^\times)^2$ has order divisible by $2^{m+2}$.
\endproclaim
\demo{Proof}
By Kummer theory, $F\bigl(\sqrt{\Gamma_F^S}\bigr)/F$ is the maximal multiquadratic extension of $F$ which is 
unramified outside of $S\cup S_2$. Hence  $F\bigl(\sqrt{\Gamma_F^S}\bigr)$ contains $\E(\sqrt{-1})$. Kummer theory
also gives the degree of the former field over $F$ as
$|\Gamma_F^S/(F^\times)^2|$, while the latter has degree $2^{m+1}$ over $F$. The result follows.
\qed
\enddemo

Combining Proposition 6.1 and Lemma 6.2 yields the following.

 \proclaim{Corollary 6.3}
The group $_2K_2(\Oc_F^S)$ has order divisible by $2^{m}$.
If $\E(\sqrt{-1})$ is not the maximal multiquadratic extension of $F$ which is unramified outside of 
$S\cup S_2$, then  $K_2(\Oc_F^S)[2]$ has order divisible by $2^{m+1}$.
\endproclaim

For each relative quadratic extension $E$ of $F$ inside $\E$, we will also need to know the order and the
$2$-rank of 
$K_2(\Oc_E^{S})^{1-\tau}$, where $\tau$ is the nontrivial automorphism of $E$ over $F$. 
These can be described in terms of the homomorphism $\iota_*: K_2(\Oc_F^S)\rightarrow K_2(\Oc_E^S)$ obtained
by functorality from the inclusion of $\Oc_F^S$ in $\Oc_E^S$.
Note that the image of $K_2(\Oc_F^S)$ under $\iota_*$ is contained in $K_2(\Oc_E^{S})^{\langle \tau \rangle}$, the subgroup of $\tau$-fixed elements of $K_2(\Oc_E^{S})$. Following standard practice, we let $\ker(\iota_*)$
and $\coker(\iota_*)$ denote
the kernel and cokernel of $\iota_*$ as a map from $K_2(\Oc_F^S)$ to $K_2(\Oc_E^{S})^{\langle \tau \rangle}$.

\proclaim{Proposition 6.4} 
$$|K_2(\Oc_E^{S})^{1-\tau}|=\dfrac{|K_2(\Oc_E^{S})|}{|K_2(\Oc_F^{S})|}\dfrac{|\ker(\iota_*)|}{|\coker(\iota_*)|}.$$
\endproclaim
\demo{Proof}
This follows directly from the exact sequences
$$ 0\rightarrow K_2(\Oc_E^S)^{\langle\tau\rangle}\rightarrow K_2(\Oc_E^S)
\overset{1-\tau}\to\rightarrow K_2(\Oc_E^S)^{1-\tau}\rightarrow 0$$ and
$$0\rightarrow \ker(\iota_*)\rightarrow K_2(\Oc_F^S)\overset{\iota_*}\to\rightarrow K_2(\Oc_E^S)^{\langle\tau\rangle}
\rightarrow \coker(\iota_*) \rightarrow 0. \qed$$
\enddemo

We pause to note the consequences of Kahn's result for our specific situation.

\proclaim{Remark 6.5}{\rm
\item{a.} In the current setting, $r_\infty(E/F)=0$
and
$K_2(\Oc_E^S)_{\langle \tau \rangle}=K_2(\Oc_E^S)/K_2(\Oc_E^S)^{1-\tau}$. 
Corollary 3.6 then gives $K_2(\Oc_E^S)/K_2(\Oc_E^S)^{1-\tau}
\cong K_2(\Oc_F^S)$, with the isomorphism induced by the transfer map
$\Tr$. Thus $|K_2(\Oc_E^S)^{1-\tau}|=\dfrac{|K_2(\Oc_E^S)|}{|K_2(\Oc_F^S)|}$.
\smallskip
\item{b.} Comparing this equation with that of Proposition 6.4 shows that
$|\ker(\iota_*)|=|\coker(\iota_*)|$. This also follows from the Hochschild-Serre
spectral sequence. Interestingly, we do not seem to make direct use of this equality.
}
\endproclaim

\proclaim{Proposition 6.6}
$$\rank_2(K_2(\Oc_E^S)^{1-\tau})\geq \rank_2(K_2(\Oc_F^S))-\rank_2(\coker(\iota_*)),$$
and $\coker(\iota_*)$ has exponent 2.
\endproclaim
\demo{Proof}
$$K_2(\Oc_E^{S})^{1-\tau}\cong K_2(\Oc_E^{S})/K_2(\Oc_E^{S})^{\langle \tau \rangle}
\cong \bigl(K_2(\Oc_E^{S})/K_2(\Oc_E^{S})^{1+\tau}\bigr)/
\bigl(K_2(\Oc_E^{S})^{\langle \tau \rangle}/K_2(\Oc_E^{S})^{1+\tau}\bigr)$$
so
$\rank_2(K_2(\Oc_E^S)^{1-\tau}\geq 
\rank_2\bigl(K_2(\Oc_E^{S})/K_2(\Oc_E^{S})^{1+\tau}\bigr)
-\rank_2\bigl(K_2(\Oc_E^{S})^{\langle \tau \rangle}/K_2(\Oc_E^{S})^{1+\tau}\bigr)$.
Now $K_2(\Oc_E^{S})/K_2(\Oc_E^{S})^{1+\tau}$ modulo squares is identical to
$K_2(\Oc_E^{S})/K_2(\Oc_E^{S})^{1-\tau}$ modulo squares. So the 2-rank of the former
equals the 2-rank of the latter, which is $\rank_2(K_2(\Oc_F^S))$, by the isomorphism in Remark 6.5a.
For the quotient
$K_2(\Oc_E^{S})^{\langle \tau \rangle}/K_2(\Oc_E^{S})^{1+\tau}$,
note first that it has exponent 2, as 2 acts on it as $1+\tau$. 
Now $1+\tau$ factors as $\iota_*\circ \Tr$ in this setting, and $\Tr$ defines
a surjection onto $K_2(\Oc_F^S)$, again by the use of Kahn's theorem as spelled out in Remark 6.5.
Hence this quotient is exactly $\coker(\iota_*)$.
\qed
\enddemo

\heading{VII. Proof of the Annihilation Theorem}
\endheading

When assuming the truth of the Birch-Tate conjecture for $F$ and the intermediate relative quadratic
extensions $E$ in $\E$, put $\delta=0$ if $\E(\sqrt{-1})$ is not the
maximal multiquadratic extension of $F$ that is unramified outside of 
$S\cup S_2$; and put $\delta=1$ when it is. When not assuming the truth of the Birch-Tate conjecture for these fields,
we define $\delta$ so that the order of the 2-Sylow subgroup of
$K_2(\Oc_\E^S)$ is $2^\delta$. 
 
We begin the proof of Theorem 4.1 by showing that the first term of 
$(-1)^{|S|}2^\delta\alpha_{\E/F}^S(\q,1)$ from Proposition 5.1, namely the term
$\dfrac{(\N\q^2-1)}{w_2(F)}\cdot \dfrac{2^\delta|K_2(\Oc_F^S)|}{2^m}T_{\E/F}$
corresponding to the principal character, annihilates 
$K_2(\Oc_\E^S)$. First, $\dfrac{\N\q^2-1}{w_2(F)}\in \Z$, by Proposition 1.1.
 Next, $T_{\E/F}$ factors as the composition of  $\Tr=\Tr_{\E/F}: K_2(\Oc_\E^S) \rightarrow K_2(\Oc_F^S)$     
and $\iota_* : K_2(\Oc_F^S)\rightarrow K_2(\Oc_\E^S)$ obtained by functorality from the inclusion of $F$ in $E$. 
The fact from Corollary 6.3 that $\rank_2(K_2(\Oc_F^S))\ge m+1-\delta$
implies that $\dfrac{2^\delta|K_2(\Oc_F^S)|}{2^m}$ is an integer which annihilates $K_2(\Oc_F^S)$; it thus annihilates the homomorphic
image $\iota_*(K_2(\Oc_F^S))$, and so also annihilates $T_{\E/F}K_2(\Oc_\E^S)\subset \iota_*(K_2(\Oc_F^S)$.
Thus $\dfrac{2^\delta|K_2(\Oc_F^S)|}{2^m}T_{\E/F}$ annihilates $K_2(\Oc_\E^S)$, and hence so does the integer multiple 
which is the first term of $(-1)^{|S|}2^\delta\alpha_{\E/F}^S(\q,1)$.

 For the other terms, we need some additional integrality results.

\proclaim{Lemma 7.1} Suppose that $E/F$ is a relative quadratic extension of totally real fields and that
$\q$ is a prime of $F$ which does not divide $w_2(E)$ and is unramified in
$E/F$. 
\item{a.} If $\q$ splits in $E/F$, then $w_2(E)|(\N\q^2-1)$.
\item{b.} If $\q$ remains inert in $E/F$, then
$w_2(E)|(\N\q^2+1)w_2(F)$.
\endproclaim
\demo{Proof} \item{a.} If $\q$ splits in E and $\fQ$ lies over $\q$,
then $\N\q=\N\fQ$ and the result follows from Proposition 1.1 (applied to the field $E$).  
\item{b.} Note that $E(W_2(E))/F$ is the composite of the abelian
extensions $E/F$ and $F(W_2(E))/F$, and so has an abelian Galois group $A$. 
The subgroup $H=\Gal(E(W_2(E))/E$ has exponent 2 by the definition of $W_2(E)$, and is of index 2 in $A$. 
Let $\sigma_{\q}\in A$ be the Frobenius at $\q$ so that its restriction to $E$ generates $\Gal(E/F)$. It
follows that $A^2=\langle {\sigma_\q}^2\rangle$, so that
$W_2(E)^{\sigma_\q^2+1}$ is fixed by $A^2$ and hence consists of
roots of unity satisfying the definition of $W_2(F)$. But $W_2(E)^{\sigma_\q^2+1}=W_2(E)^{\N\q^2+1}$ so that  
$W_2(E)^{\N\q^2+1}\subset W_2(F)$ and consequently $(\N\q^2+1)w_2(F)$
annihilates $W_2(E)$. The result follows from the cyclicity of $W_2(E)$.
\qed
\enddemo

We now show that at the prime 2, a stronger statement holds. Let 
 $\Q_r$ denote the $r$th layer of the cyclotomic
$\Z_2$ extension of $\Q$; it may be defined as $\Q_r=\Q(\mu_{2^{r+2}})^+$, the maximal real subfield of the
field of $2^{r+2}$th roots of unity. Thus $\Q_0=\Q$ and $\Q_1=\Q(\sqrt{2})$. Also let
$\Q_\infty$ denote the union of the $\Q_r$ for all natural numbers $r$. 

\proclaim{Lemma 7.2} If $M$ is a real field, and $M\cap \Q_\infty = \Q_r$, then
$2^{r+3}$ exactly divides $w_2(M)$.
\endproclaim
\demo{Proof}
The extension $M(\mu_{2^{r+3}})/M$ is the composite of the two extensions $M(\sqrt{-1})$ and 
$M\cdot \Q_{r+1}$, both of which are relative quadratic. Hence this extension is biquadratic, and
$\mu_{2^{r+3}}\subset W_2(M)$. Thus $2^{r+3} | w_2(M)$. On the other hand,
the extension $\Gal(M(\mu_{2^{r+4}})/M)\cong \Gal(\Q(\mu_{2^{r+4}})/\Q_r)\supset
\Gal(\Q(\mu_{2^{r+4}})/\Q(\mu_{2^{r+2}}))$, which is cyclic of order 4.
So $\mu_{2^{r+4}}\not\subset W_2(F)$ and $2^{r+4} \nmid w_2(F)$.
\qed
\enddemo

\proclaim{Lemma 7.3} 
The kernel of the map $\iota_* : K_2(\Oc_F^S)\rightarrow K_2(\Oc_E^{S})$
has order 2 unless $E=F_1$, the first layer of the cyclotomic $\Z_2$-extension of $F$,
in which case it has order 1.
\endproclaim
\demo{Proof}
The composition $\Tr\circ \iota_*$ multiplies each element by the degree $[E:F]=2$, so 
$\ker(\iota_*)\subset K_2(\Oc_F^S)[2]$. The proof of Proposition 6.1 shows that each
element of $K_2(\Oc_F^S)[2]$ is represented in the form $\{-1,a\}_F$ for some 
$a \in\Gamma_F^S$ (as defined in section 6).
The image of this element under $\iota_*$ is represented by $\{-1,a\}_E$. By definition,
$\{-1,a\}_E$ is trivial if and only if $a\in D_E$. As seen in Proposition 6.1, 
this occurs if and only if 
$E(\sqrt{a})$ is contained in the first layer $E_1$ of the cyclotomic
$\Z_2$-extension of $E$, and this is clearly equivalent to the condition that
$F(\sqrt{a})\subset E_1$. Since $E_1/E$ is unramified outside of $S_2$, we conclude 
that for $a \in F^\times$,
the condition $\{-1,a\}_F\in \ker(\iota_*)$
is equivalent to $F(\sqrt{a})\subset E_1$. Let $\A$ denote the set of $a\in F^\times$
satisfying this condition.

First suppose that $E\neq F_1$. Then $E_1=E\cdot F_1$ is a biquadratic extension of $F$ which is unramified outside
of $S\cup S_2$. In this case, $\{-1,a\}\in \ker(\iota_*)$ if and only if $F(\sqrt{\a})$ lies in this biquadratic extension.
Kummer theory then gives us that $\A/(F^\times)^2$ has order 4, and so by Proposition 6.1, 
$\ker(\iota_*)\cong \A/D_F \cong (\A/(F^\times)^2)/(D_F/(F^\times)^2)$ has order 2.
  
Now suppose that $E=F_1$. Then $E_1=F_2$ is a cyclic extension
of degree 4 over $F$, and the only quadratic extension it contains is $E=F_1$.
In this case we see that $\{-1,a\}_F\in \ker(\iota_*)$
if and only if
$F(\sqrt{a})\subset F_1$. But, as mentioned in Proposition 6.1, this occurs if and only if $a\in D_F$, and then
$\{-1,a\}_F$ is trivial. Thus $\ker(\iota_*)$ is trivial in this case.
\qed
\enddemo

Let $\psi$ denote the non-trivial character of $\Gal(E/F)=\langle\tau\rangle$.

\proclaim{Proposition 7.4} The quantity  $\dfrac{(\N\q^2-\psi(\sigma_\q))w_2(F)}{|\ker(\iota_*)|w_2(E)}$
is integral.
\endproclaim
\demo{Proof}
If $\sigma_\q\in\Gal(E/F)$ is trivial, then 
$\dfrac{(\N\q^2-\psi(\sigma_\q))}{w_2(E)}\in \Z$ by Lemma 7.1a. Furthermore,
$\dfrac{w_2(F)}{|\ker(\iota_*)|}\in \Z$ by Lemma 7.3, since $2|w_2(F)$ (indeed $2^3|w_2(F)$).
If $\sigma_\q=\tau$, then 
$\dfrac{(\N\q^2-\psi(\sigma_\q))w_2(F)}{w_2(E)}$ is integral by Lemma 7.1b, so
we need only check the integrality of 
$\dfrac{(\N\q^2-\psi(\sigma_\q))w_2(F)}{|\ker(\iota_*)|w_2(E)}$
at 2. If $\dfrac{w_2(F)}{w_2(E)}$ is integral at 2 then we're done 
by Lemma 7.3.
If $\dfrac{w_2(F)}{w_2(E)}$ is not integral at 2 then $\Q_t=E\cap\Q_\infty \neq F\cap\Q_\infty=\Q_r$, by 
Lemma 7.2. Then $E=F\cdot\Q_t$ and $\Q_r=F\cap\Q_t$ so $[\Q_t:\Q_r]=[E:F]=2$.
Thus $t=r+1$ and $E=F_1$, the first layer of the cyclotomic $\Z_2$ extension of $F$.
Now by Lemma 7.3, $|\ker(\iota_*)|=1$; and by Lemma 7.2 again, $\dfrac{2w_2(F)}{w_2(E)}$ is integral at 2
and hence so is $\dfrac{(\N\q^2-\psi(\sigma_\q))w_2(F)}{w_2(E)}$ because
$\N\q^2-\psi(\sigma_\q)$ is even. 
\qed
\enddemo

 Now we are ready to complete the proof of 
Theorem 4.1 by showing that each of the remaining terms
$\pm 2^\delta \dfrac{w_2(F)}{w_2(E)}\dfrac{|K_2(\Oc^S_{E})|}
{2^m|K_2(\Oc^S_F)|}(\N\q^2-\sigma_\q)(1-\tau) T_{\E/E}$
of $2^\delta\alpha_{\E/F}^S(\q,1)$ annihilates $K_2(\Oc_\E^S)$.

Since $T_{\E/E}: K_2(\Oc_\E^S)\rightarrow K_2(\Oc_\E^S)$ factors 
as a $G$-module homomorphism through
$K_2(\Oc_E^S)$, it will suffice to show that  
$2^\delta \dfrac{w_2(F)}{w_2(E)}\dfrac{|K_2(\Oc^S_{E})|}
{2^m|K_2(\Oc^S_F)|}(\N\q^2-\sigma_\q)(1-\tau)$
annihilates $K_2(\Oc_E^S)$. This depends only on $\sigma_\q$ as an automorphism of $E$, for
which $\sigma_\q (1-\tau)=\psi(\sigma_\q) (1-\tau)$. We are then reduced to showing that
$$2^\delta \dfrac{w_2(F)}{w_2(E)}\dfrac{|K_2(\Oc^S_{E})|}
{2^m|K_2(\Oc^S_F)|}(\N\q^2-\psi(\sigma_\q))(1-\tau)$$
annihilates
$K_2(\Oc_E^S)$, or equivalently that 
$2^\delta \dfrac{w_2(F)}{w_2(E)}\dfrac{|K_2(\Oc^S_{E})|}
{2^m|K_2(\Oc^S_F)|}(\N\q^2-\psi(\sigma_\q))$ annihilates
$K_2(\Oc_E^S)^{(1-\tau)}$.
By Proposition 7.4, $\dfrac{(\N\q^2-\psi(\sigma_\q))w_2(F)}{|\ker(\iota_*)|w_2(E)}$
is integral, leaving us to show that 
$$2^\delta |\ker(\iota_*)|\dfrac{|K_2(\Oc^S_{E})|}
{2^m|K_2(\Oc^S_F)|}$$ annihilates
$K_2(\Oc_E^S)^{(1-\tau)}$. 
At this point, Proposition 6.4 states that  $K_2(\Oc_E^S)^{(1-\tau)}$
has order equal to $\dfrac{|K_2(\Oc_E^{S})|}{|K_2(\Oc_F^{S})|}\dfrac{|\ker(\iota_*)|}{|\coker(\iota_*)|}$,
while Proposition 6.6 and Corollary 6.3 give the result that $K_2(\Oc_E^S)^{(1-\tau)}$
has 2-rank greater than or equal to 
$$\rank_2(K_2(\Oc_F^S))-\rank_2(\coker(\iota_*))\geq
m+1-\delta-\rank_2(\coker(\iota_*)),$$ 
and $|\coker(\iota_*)|=2^{\rank_2(\coker(\iota_*))}$.  
Thus $K_2(\Oc_E^S)^{(1-\tau)}$ is annihilated by
$$\dfrac{|K_2(\Oc_E^{S})|}{|K_2(\Oc_F^{S})|}\dfrac{|\ker(\iota_*)|}{|\coker(\iota_*)|}
\dfrac{|\coker(\iota_*)|}{2^{m-\delta}}=
\dfrac{2^\delta |K_2(\Oc_E^S)|}{|K_2(\Oc_F^S)|}\dfrac{|\ker(\iota_*)|}{2^m},$$
as desired.

\heading{ACKNOWLEDGEMENTS}
\endheading
We thank Manfred Kolster for sharing his wealth of knowledge on the subject of $K$-theory and arithmetic.

\enddocument
\end